\documentclass{amsart}

\theoremstyle{plain}
\newtheorem{prop}{Proposition}[section]

\theoremstyle{definition}
\newtheorem{defn}[prop]{Definition}
\newtheorem{exmp}[prop]{Example}
\newtheorem*{prob}{Problem}

\theoremstyle{remark}
\newtheorem*{note}{Note}

\begin{document}

\title[Examples of $k$-nets]{Old and new examples of $k$-nets in $\mathbb{P}^2$}

\author{Janis Stipins}

\address{Department of Mathematics
\\University of Michigan
\\525 East University Avenue
\\Ann Arbor, MI 48109-1109}

\email{jstipins@umich.edu}

\dedicatory{Dedicated to Professor Sergey Yuzvinsky on the occasion of his seventieth birthday.}


\begin{abstract}
In this paper, we present a number of examples of $k$-nets, which are
special configurations of lines and points in the projective plane.
Such a configuration can be regarded as the union of $k$ 
completely reducible 
elements of a pencil of complex plane curves; equivalently it can be
regarded as a set of $k$ polygons in the complex projective
plane that satisfy a condition of mutual perspectivity and nondegenerate
intersection.  
For each example, we describe its construction, combinatorial properties,
and parameter space.
Most of the examples are historical, although perhaps not very
well-known; our only essentially new example is a $3$-net of pentagons which
does not realize a group.  The existence of this example settles
a question posed by S.\ Yuzvinsky \cite{yuz:nets}.
\end{abstract}

\maketitle

\section{Definition and basic terminology}

Let us begin by recalling the combinatorial definition of a $k$-net
in a projective plane:

\begin{defn}  
Let $k$ be a positive integer and let $P$ be a projective plane.  
A $k$-net in $P$ is a $(k+1)$-tuple 
$(\mathcal{A}_1, \ldots, \mathcal{A}_k, \mathcal{X})$, where each
$\mathcal{A}_i$ is a nonempty finite set of lines of $P$ and $\mathcal{X}$ is
a finite set of points of $P$, satisfying the following conditions:
\begin{enumerate}

\item The $\mathcal{A}_i$ are pairwise disjoint.

\item If $i \neq j$, then the intersection point of any line in 
$\mathcal{A}_i$ with any line in $\mathcal{A}_j$ belongs to $\mathcal{X}$.

\item Through every point in $\mathcal{X}$ there passes
exactly one line of each $\mathcal{A}_i$.

\end{enumerate}
\end{defn}

It is easy to see that this definition is uninteresting when $k < 3$.
For example, any two disjoint nonempty finite sets of lines $\mathcal{A}_1$ and
$\mathcal{A}_2$ form a $2$-net, with 
$\mathcal{X} = \left\{ \, l \cap l' \, | \, l \in \mathcal{A}_1 , 
\, l' \in \mathcal{A}_2 \right\}$.  When $k \ge 3$, however, there
are very strong restrictions on the combinatorics of a $k$-net:

\begin{prop}
Let $(\mathcal{A}_1, \ldots, \mathcal{A}_k, \mathcal{X})$ 
be a $k$-net, with $k \ge 3$.
Then every $\mathcal{A}_i$ has the same cardinality.
Furthermore, the points of $\mathcal{X}$ are exactly the 
intersections of lines of $\mathcal{A}_i$ with lines of
$\mathcal{A}_j$, for any $i \neq j$.  Thus 
$| \mathcal{X} | = | \mathcal{A}_1 |^2 $.
\end{prop}
\begin{proof}
For every $i < j$, let $n_{ij}$ 
be the number of intersections formed by lines of
$\mathcal{A}_i$ with lines of $\mathcal{A}_j$.
Since every such intersection
must be a point of $\mathcal{X}$, we have
\[
	| \mathcal{X} | \ge \mathrm{max} \left\{ n_{ij} \right\} .
\]
On the other hand, for every $i < j$, 
each point of $\mathcal{X}$ must be an intersection
point of a line of $\mathcal{A}_i$ with a line of $\mathcal{A}_j$.  
So we also have
\[
	\mathrm{min} \left\{ n_{ij} \right\} \ge | \mathcal{X} |.
\]
We conclude that every $n_{ij}$ is equal to $| \mathcal{X} |$.  
However, the condition
that only one line of each $\mathcal{A}_i$ passes through each point
of $\mathcal{X}$ implies that
\[
	n_{ij} = | \mathcal{A}_i | \cdot | \mathcal{A}_j | .
\]
Since $k \ge 3$, these products are all equal only if every
$\mathcal{A}_i$ has the same cardinality.  The rest of the claim
follows easily.
\end{proof}

\begin{note}
Because of this result, from now on we will always assume 
$k \ge 3$ in a $k$-net.
\end{note}

Yuzvinsky \cite{yuz:nets} denotes the common value of $| \mathcal{A}_i |$
by $m$, and calls this number the \emph{order} of the $k$-net.  This is
sensible terminology, because we will see shortly that $m$ is
the common order of a set of Latin squares associated to the $k$-net.
However, we will
tend to regard the lines in a given $\mathcal{A}_i$ as the components of a
completely reducible plane curve, and so we will use slightly different
terminology: we will denote the common value of $| \mathcal{A}_i |$ by
$d$, and call it the \emph{degree} of the $k$-net.  We will continue the
analogy with plane curves, and refer to a $k$-net of degree three as a
``$k$-net of cubics,'' etc.

Following 
Yuzvinsky, we will use the notation $(k,d)$-net for a $k$-net of 
degree $d$.  Also, we will only consider $(k,d)$-nets with $d > 1$.

\section{$k$-nets and Latin squares}

Let us recall a definition from elementary projective geometry:

\begin{defn}
Let $l_1, \ldots, l_d$ and $l'_1, \ldots, l'_d$ be two disjoint sets
of lines in a projective plane.  We say that the $d$-gons formed 
by $\left\{ l_i \right\}$ and $\left\{ l'_i \right\}$ are 
\emph{perspective from the line $m$} when there exists a permutation
$\sigma \in S_d$ such that the intersection points
$l_i \cap l'_{\sigma(i)}, 1 \le i \le d$, all lie on $m$.
\end{defn}

\begin{note}
Here, and throughout, we use the term \emph{$d$-gon} to mean any 
union of $d$ lines in the plane; there is no assumption made about
nondegeneracy.  The $d$ lines could be concurrent, in general position,
or anything in between.
\end{note}

Now we can state the following proposition, whose proof is immediate:

\begin{prop}\label{prop:3net-pers}
In a $(3,d)$-net $(\mathcal{A}_1, \mathcal{A}_2, \mathcal{A}_3, \mathcal{X})$,
any two of the $d$-gons formed by the $\mathcal{A}_i$ are perspective from
every side of the third.
\end{prop}

There is a convenient way to record the perspectivities in a $(3,d)$-net.
Choose a labeling for the lines of $\mathcal{A}_1$ and $\mathcal{A}_2$:
\begin{eqnarray*}
	\mathcal{A}_1 = \{ l_{11}, \ldots, l_{1d} \} \\
	\mathcal{A}_2 = \{ l_{21}, \ldots, l_{2d} \}
\end{eqnarray*}
Now we know that the line $l_{11}$ establishes a perspectivity between
the $d$-gons formed by $\mathcal{A}_2$ and $\mathcal{A}_3$.  We label
the lines of $\mathcal{A}_3$ so that the intersections
\[
l_{21} \cap l_{31},\,\, l_{22} \cap l_{32}, \,\,\ldots, \,\,l_{2d} \cap l_{3d}
\]
all lie on $l_{11}$.  Now we can create a $d \times d$ matrix $M_3$
(notation to be explained shortly)
such that $(M_3)_{ij} = n$ if and only if $l_{1i}, l_{2j}$, and $l_{3n}$
are coincident.  

\begin{prop}\label{prop:matrix-latin}
The matrix $M_3$ is a Latin square of order $d$.  The $i$-th row of
$M_3$ describes the permutation corresponding to the perspectivity
that $l_{1i}$ establishes between the $d$-gons corresponding to
$\mathcal{A}_2$ and $\mathcal{A}_3$, with labels chosen as described.
\end{prop}

A different choice of labels for any of the
$\mathcal{A}_i$ leads to a Latin square which is said to be
\emph{isotopic} to $M_3$.  Every Latin square is the multiplication
table of a quasigroup; our method of labeling ensures that the
first row of $M_3$ is just $(1,2,\ldots,d)$, in which case the
quasigroup is in fact a \emph{loop} \cite{den:latin}.

To use Latin squares in describing $k$-nets when $k > 3$, we need 
one more elementary concept:

\begin{defn}
Let $L$ and $M$ be Latin squares of order $d$.
$L$ and $M$ are an \emph{orthogonal pair} if every element
of $\{1,2,\ldots,d\} \times \{1,2,\ldots,d\}$ appears exactly once
in the set
$\left\{ ( L_{ij}, M_{ij} ) \right\}$.
\end{defn}
\begin{exmp}\label{ex:43net}
The following two Latin squares form an orthogonal pair:
\[
\left[\begin{array}{ccc}
1 & 2 & 3 \\
2 & 3 & 1 \\
3 & 1 & 2
\end{array}\right] \,\,\,
\left[\begin{array}{ccc}
1 & 2 & 3 \\
3 & 1 & 2 \\
2 & 3 & 1
\end{array}\right]
\]
\end{exmp}

\begin{defn}
A set of Latin squares of order $d$ forms an \emph{orthogonal set}
if every pair of Latin squares in the set is orthogonal.
\end{defn}

\begin{exmp}\label{ex:54net}
The following three Latin squares form an orthogonal set:
\[
\left[\begin{array}{cccc}
1 & 2 & 3 & 4 \\
2 & 1 & 4 & 3 \\
3 & 4 & 1 & 2 \\
4 & 3 & 2 & 1
\end{array}\right] \,\,\,
\left[\begin{array}{cccc}
1 & 2 & 3 & 4 \\
3 & 4 & 1 & 2 \\
4 & 3 & 2 & 1 \\
2 & 1 & 4 & 3
\end{array}\right] \,\,\,
\left[\begin{array}{cccc}
1 & 2 & 3 & 4 \\
4 & 3 & 2 & 1 \\
2 & 1 & 4 & 3 \\
3 & 4 & 1 & 2
\end{array}\right]
\]
\end{exmp}

For much more information about orthogonal sets of Latin squares, 
see \cite{den:latin}.  For our purposes, orthogonal sets of Latin 
squares allow us to describe the combinatorial properties of $k$-nets.
The following proposition is a generalization of 
Proposition~\ref{prop:3net-pers}, and hopefully explains the
notation used in Proposition~\ref{prop:matrix-latin}:

\begin{prop}\label{prop:knet-pers}
In $(k,d)$-net $(\mathcal{A}_1, \ldots, \mathcal{A}_k, \mathcal{X})$,
any two $d$-gons formed by the $\mathcal{A}_i$ are perspective from
every side of every other $d$-gon.
For each $d$-gon $\mathcal{A}_i$, it is
possible to choose a labeling $\mathcal{A}_i = \{ l_{i1}, \ldots, l_{id} \}$
and to construct $k-2$ matrices $M_m, 3 \le m \le k,$ such that:
\begin{enumerate}
\item The first row of every $M_m$ is $(1,2,\ldots,d)$.

\item $(M_m)_{ij} = n$ if and only if the lines
$l_{1i}, l_{2j}$, and $l_{mn}$ are coincident.

\item The collection $M_3, M_4, \ldots, M_k$ is an orthogonal set of 
Latin squares of order~$d$.
\end{enumerate}  
In other words, the matrix $M_m$ describes the perspectivities between
the $d$-gons corresponding to 
$\mathcal{A}_1$, $\mathcal{A}_2$, and $\mathcal{A}_m$, with labels
chosen as described.
\end{prop}

There is a sort of converse to the preceding proposition, which we will
use later as a recipe for constructing examples of $k$-nets.  We start
with a simple observation:

\begin{prop}
Let $\mathcal{X}$ be a set of $d^2$ distinct points in a projective plane.  A
necessary and sufficient condition for $\mathcal{X}$ to be the set of
points of a $(k,d)$-net is that $\mathcal{X}$ can be partitioned $k$
different ways into $d$ sets of $d$ collinear points.  Here each 
$\mathcal{A}_i$
is the set of $d$ lines, each containing $d$ points, corresponding to the
$i$-th such partition.
\end{prop}

Now it is easy to see the following converse of 
Proposition~\ref{prop:knet-pers}:

\begin{prop}\label{prop:knet-construct}
Suppose $\mathcal{A}_1$ and $\mathcal{A}_2$ are the sets of lines of two
$d$-gons which intersect in $d^2$ distinct points; denote by $\mathcal{X}$
this set of intersection points.  Choose any labeling of the lines of 
$\mathcal{A}_1$ and $\mathcal{A}_2$:
\[
\mathcal{A}_1 = \{ l_{11}, \ldots, l_{1d} \}, \,\,
\mathcal{A}_2 = \{ l_{21}, \ldots, l_{2d} \} .
\]
Suppose $\mathcal{A}_3, \ldots, \mathcal{A}_k$ 
are sets of lines which establish perspectivities
between the lines of $\mathcal{A}_1$ and $\mathcal{A}_2$ according to an
orthogonal set of Latin squares
Latin square $M_3, \ldots, M_k$, 
as described in Proposition~\ref{prop:knet-pers}.
Then $(\mathcal{A}_1, \ldots, \mathcal{A}_k, \mathcal{X})$ is a 
$(k,d)$-net.
\end{prop}
\begin{proof}
The only observation to make here is that the perspectivities described
by each $M_m$ correspond to a partition of $\mathcal{X}$ as described
in the preceding proposition; the fact that each $M_m$ is a Latin square
is exactly the condition that every point of $\mathcal{X}$ is contained
in exactly one line of each $\mathcal{A}_m$.
\end{proof}

\begin{exmp}
One application of the ideas of Propositions~\ref{prop:knet-pers} and
\ref{prop:knet-construct} is the following fact: There is a correspondence
between finite projective planes of order~$n$ and $(n+1,n)$-nets.
To see this, let $P$ be a projective plane of order~$n$, and
choose any line $l$ of $P$; it
contains $n+1$ points, $X_1, \ldots, X_{(n+1)}$.  
Let $\mathcal{A}_i$ be the set of $n$ lines through $X_i$ other than $l$,
and let $\mathcal{X}$ be all the points of $P$ not lying on $l$.  Then
$(\mathcal{A}_1, \ldots, \mathcal{A}_{(n+1)}, \mathcal{X})$ is a 
$(n+1, n)$-net.  This example is maximal in the sense that $n-1$ is the
largest possible cardinality of an orthogonal set of Latin squares of
order~$n$.  It turns out, conversely,
that any orthogonal set of $n-1$ Latin 
squares of order $n$ can be used to construct a finite projective plane
of order~$n$; the idea of the construction is more or less the same as 
the idea of Proposition~\ref{prop:knet-construct}.  
In fact, the Latin squares in 
Example~\ref{ex:43net} describe a $(4,3)$-net which can be used to 
construct a projective plane of order three; and
the Latin squares in Example~\ref{ex:54net} describe a $(5,4)$-net
which can be used to construct a projective plane of order four.
See \cite{den:latin} for details and proofs of these facts.
\end{exmp}

Like Yuzvinsky, we are most interested in $k$-nets which can be 
realized in the complex projective plane.  In general, the
existence of an orthogonal set of $k-2$ Latin squares of order~$d$
is not sufficient for the realizability of a $(k,d)$-net in the complex
projective plane; for example, in the next section we will see that
a $(5,4)$-net corresponding to Example~\ref{ex:54net} is not so 
realizable, for topological reasons.  However,
Proposition~\ref{prop:knet-pers} shows that
\emph{nonexistence} of an orthogonal set of $k-2$ Latin squares 
of order~$d$ is certainly enough to prove that a $(k,d)$-net cannot
exist.  For example, an old problem of Euler (\cite{den:latin}, p.\ 156)
was solved by proving that there is
no pair of orthogonal Latin squares of order six; this means that we
cannot hope to find a $(4,6)$-net in the complex projective plane, or
in any other projective plane.

\section{$k$-nets and pencils of curves in $\mathbb{P}^2(\mathbb{C})$}

Now let us restrict our attention to projective planes of the 
form $\mathbb{P}^2(K)$, where $K$ is an algebraically closed field
(not necessarily of characteristic zero).  In this case, Max Noether's
Fundamental Theorem \cite{fulton:curves} gives the following result:

\begin{prop}
Let $C_1, \ldots, C_k$ be a collection of reduced, completely reducible
curves of degree $d$ in $\mathbb{P}^2(K)$.  Let $\mathcal{A}_i$ be the
union of the lines forming the curve $C_i$.  Then the $\mathcal{A}_i$
form a $(k,d)$-net if and only if the $C_i$ belong to a single pencil of
plane curves.
\end{prop}

This means that if we are searching for examples of $(k,d)$-nets, we
can start by analyzing pencils of degree~$d$ curves with $d^2$ distinct
base points.
When the characteristic of $K$ is zero, we can say a little more about
these pencils:

\begin{prop}
Suppose $\mathrm{char}(K) = 0$.  Then the general member of the 
pencil of curves in the preceding proposition is smooth and reduced.
Conversely, suppose we are given a pencil of plane curves of 
degree $d$, whose general member is smooth and reduced; then the completely
reducible members of this pencil must also be reduced, and any $k$
such members form a $(k,d)$-net.
\end{prop}

For the remainder of this paper, we will
restrict our attention to the very interesting case
that $K = \mathbb{C}$.  In this case, we can use topological arguments
which are not available over an arbitrary field.  To wit: 
Given a $(k,d)$-net in 
$\mathbb{P}^2 = \mathbb{P}^2(\mathbb{C})$, the preceding propositions
enable us to recognize the $d^2$ points of $\mathcal{X}$ as the
base points of a pencil of degree $d$ curves, whose general element
is smooth.  Blowing up these base points gives a rational surface
fibered over $\mathbb{P}^1$, and a straightforward argument using
Euler characteristics puts very strong restrictions on the possibilities
for $k$ and $d$.  Yuzvinsky \cite{yuz:nets} carries out this argument 
to conclude the following:

\begin{prop}
For an arbitrary $(k,d)$-net in $\mathbb{P}^2$, the only possible
values for $(k,d)$ are: $(k=3, d\ge 2), (k=4, d\ge 3), (k=5, d\ge 6)$.
\end{prop}

In fact, it is conjectured that the only $4$-net in $\mathbb{P}^2$
is the $(4,3)$-net described below, and it is also conjectured that 
there are no $5$-nets in $\mathbb{P}^2$ whatsoever.  A detailed 
study of these conjectures is the topic of a future paper; for now,
we move on to particular examples of $k$-nets which are realizable
in the complex projective plane.

\section{Examples}

\subsection{$3$-nets of conics}
This is an easy application of the ideas in the preceding section:
a $(3,2)$-net is the union of three reducible members of a pencil of
conics.  However, \emph{every} pencil of conics in $\mathbb{P}^2$
with four distinct base points has exactly three reducible members; if the
base points are labeled $X_1, X_2, X_3, X_4$, the three reducible members
are as follows:
\begin{eqnarray*}
	C_1 &=& \langle X_1, X_2 \rangle \cup \langle X_3, X_4 \rangle \\
	C_2 &=& \langle X_1, X_3 \rangle \cup \langle X_2, X_4 \rangle \\
	C_3 &=& \langle X_1, X_4 \rangle \cup \langle X_2, X_3 \rangle 
\end{eqnarray*}
It is easy to see that the points $X_i$ must be in general position,
and so up to projective equivalence they can be assumed to have the
following coordinates:
\[
	X_1 = \left[ 1 : 0 : 0 \right], \,\,
	X_2 = \left[ 0 : 1 : 0 \right], \,\,
	X_3 = \left[ 0 : 0 : 1 \right], \,\,
	X_4 = \left[ 1 : 1 : 1 \right].
\]
This implies that any pencil of conics with four distinct base points
is projectively equivalent to the pencil
\[
	\lambda z(x-y) + \mu y(z-x) = 0 . 
\]
Here the particular choices 
$\left[ \lambda : \mu \right] = \left[ 1 : 0 \right], 
\left[ 0 : 1 \right], \left[ 1 : 1 \right]$ correspond to
$C_1, C_2, C_3$ respectively.

We conclude that any $3$-net of conics is projectively equivalent to the
one formed by $C_1, C_2, C_3$ and $X_1, X_2, X_3, X_4$.  Finally,
we note that the lines of each $C_i$ may be labeled as in
Proposition~\ref{prop:knet-pers} to give the matrix
\[
	M_3 = \left[ \begin{array}{cc}
			1 & 2 \\
			2 & 1
			\end{array} \right] ,
\]
which is of course the only Latin square of order two.  Up to choice
of symbols, this is none other than the multiplication table for
$\mathbb{Z}/2\mathbb{Z}$.

\subsection{Two constructions of $(3,d)$-nets by Yuzvinsky}
For completeness, we recall here two constructions described in 
\cite{yuz:nets}.  First, there is an easy way 
to construct $(3,d)$-nets for any $d \ge 2$: Let $x,y,z$ be 
homogeneous coordinates on $\mathbb{P}^2$.  Then the three degenerate
curves
\begin{eqnarray*}
C_1 &=& V(x^d - y^d) \\
C_2 &=& V(y^d - z^d) \\
C_3 &=& V(z^d - x^d)
\end{eqnarray*}
form a $(3,d)$-net.  The lines can be labeled so that the corresponding
Latin square is the multiplication table for $\mathbb{Z}/d\mathbb{Z}$.
Note that each $C_i$ is a union of $d$ concurrent lines; we will avoid
this possibility in future examples.

The second construction from Yuzvinksy's paper gives a $3$-net 
corresponding to any group which can occur as a finite subgroup
of a smooth elliptic curve.  Let $G$ be such a group, and set
$d = |G|$.  Choose a smooth elliptic curve $E$ in the plane, choose one of 
the inflection points of $E$ as the identity element in the group
law, and then let $\{x_1, ..., x_d\}$ be a set of points of $E$
that realize the group $G$.  Choose a general line in the plane;
here, ``general'' means that the line should intersect $E$ in three
distinct points, which we may label as $\alpha$, $\beta$, and
$-\alpha -\beta$, using the group law on $E$.  Then it is easy 
to see that the three cosets $G + \alpha$, $G + \beta$, and
$G - \alpha - \beta$ form a configuration which is dual to a
$(3,d)$-net whose Latin square, up to labeling,
is the multiplication table for $G$.  

A rough count of parameters
shows that this second construction gives a three-dimensional family of
$(3,d)$-nets corresponding to $G$: the space of plane cubics is a
$\mathbb{P}^9$, the chosen line is a point in $(\mathbb{P}^2)^*$, and
then we take the quotient of this eleven-dimensional parameter space
by the action of $PGL(2)$, which is an eight-dimensional group; we
expect the result to be a three-dimensional family, as claimed.
However, we have already seen in the previous example that the family
of $(3,2)$-nets corresponding to $\mathbb{Z}/2\mathbb{Z}$ is
zero-dimensional: there is only one such net, up to projective
equivalence.  The explanation is that in the case that
$G = \mathbb{Z}/2\mathbb{Z}$, there is a three-dimensional family 
of plane cubics that all contain the same six
points realizing the cosets $G + \alpha$, $G + \beta$, and
$G - \alpha - \beta$ that we originally found on $E$.  A similar
discrepancy occurs in the next example.

\subsection{$3$-nets of cubics}
To find $(3,3)$-nets, we will make use of 
Proposition~\ref{prop:knet-construct}.  The perspectivities in a
$(3,3)$-net must be described by a single Latin square of order three;
up to isotopy (relabeling), there is only one such Latin square:
\[
M_3 = \left[ \begin{array}{ccc}
		1 & 2 & 3 \\
		2 & 3 & 1 \\
		3 & 1 & 2
		\end{array} \right] .
\]
We note that up to choice of symbols, this is the multiplication table
for $\mathbb{Z}/3\mathbb{Z}$.

Now we choose two sets $\mathcal{A}_1 = \{ l_{11}, l_{12}, l_{13} \}$ 
and $\mathcal{A}_2 = \{ l_{21}, l_{22}, l_{23} \}$ of three lines
each; for ease in calculation, 
we will assume that the triangles defined by these two sets are
nondegenerate (i.e., not three concurrent lines).  We may choose coordinates
on $(\mathbb{P}^2)^*$ and labels for the lines so that we have:
\begin{align*}
l_{11} &= \left[1:0:0\right] &
l_{21} &= \left[1:1:1\right] \\
l_{12} &= \left[0:1:0\right] &
l_{22} &= \left[s_0:s_1:s_2\right] \\
l_{13} &= \left[0:0:1\right] &
l_{23} &= \left[t_0:t_1:t_2\right]
\end{align*}
Apparently we have a four-parameter family of choices for $l_{22}$ and
$l_{23}$, but we have conditions to impose.  According to the Latin square
$M_3$, there must be a line $l_{31}$ containing the three intersection
points
\begin{eqnarray*}
l_{11}\cap l_{21} &=& \left[0 : 1 : -1\right], \\
l_{12}\cap l_{23} &=& \left[t_2 : 0 : -t_0\right], \\
l_{13}\cap l_{22} &=& \left[s_1 : -s_0 : 0\right] ,
\end{eqnarray*}
where here the coordinates are dual to our chosen coordinates for
$(\mathbb{P}^2)^*$.  The line $l_{31}$ exists if and only if the
corresponding determinant vanishes:
\[
\left| \begin{array}{ccc}
	0 & 1 & -1 \\
	t_2 & 0 & -t_0 \\
	s_1 & -s_0 & 0
	\end{array} \right| 
= -s_1 t_0 + s_0 t_2 = 0 ,
\]
which is to say if and only if $s_0/s_1 = t_0/t_2$.  When this 
condition is satisfied, in coordinates we have
\[
l_{31} = \left[ s_0:s_1:s_1 \right] = \left[ t_0 : t_2 : t_2 \right]. 
\]
Returning to $M_3$, we see that we
must also have a line $l_{32}$ containing the three intersection points
\begin{eqnarray*}
l_{11}\cap l_{22} &=& \left[0 : s_2 : -s_1\right], \\
l_{12}\cap l_{21} &=& \left[1 : 0 : -1\right], \\
l_{13}\cap l_{23} &=& \left[t_1 : -t_0 : 0\right] ,
\end{eqnarray*}
and such a line exists if and only if we have
\[
\left| \begin{array}{ccc}
	0 & s_2 & -s_1 \\
	1 & 0 & -1 \\
	t_1 & -t_0 & 0
	\end{array} \right| 
= -s_2 t_1 + s_1 t_0 = 0 ,
\]
which happens if and only if $s_1/s_2 = t_1/t_0$.  In this case, we 
will have
\[
l_{32} = \left[s_2:s_1:s_2\right] = \left[t_0:t_1:t_0\right] .
\]
It is not difficult to see that the conditions imposed by the existence
of $l_{31}$ and $l_{32}$ are independent, leaving us with a two-dimensional
family of pairs of triangles which are perspective from two lines.  Now we
impose the final condition: there must be a line $l_{33}$ which contains
\begin{eqnarray*}
l_{11}\cap l_{23} &=& \left[0 : t_2 : -t_1\right], \\
l_{12}\cap l_{22} &=& \left[s_2 : 0 : -s_0\right], \\
l_{13}\cap l_{21} &=& \left[1 : -1 : 0\right] .
\end{eqnarray*}
This line exists if and only if we have
\[
\left| \begin{array}{ccc}
	0 & t_2 & -t_1 \\
	s_2 & 0 & -s_0 \\
	1 & -1 & 0
	\end{array} \right| 
= -s_0 t_2 + s_2 t_1 = 0 ,
\]
which happens if and only if $s_0/s_2 = t_1/t_2$.  But this is not an
independent condition: the existence of $l_{31}$ and $l_{32}$ implies that
\[
s_0/s_2 = (s_0/s_1)(s_1/s_2) = (t_0/t_2)(t_1/t_0) = t_1/t_2 .
\]
We may use these relations to eliminate $s_2$ and $t_2$, and write 
down our family of $(3,3)$-nets in coordinates:
\begin{align*}
l_{11} &= \left[1:0:0\right] &
l_{21} &= \left[1:1:1\right] &
l_{31} &= \left[s_0:s_1:s_1\right] \\
l_{12} &= \left[0:1:0\right] &
l_{22} &= \left[s_0 t_1:s_1 t_1:s_1 t_0\right] &
l_{32} &= \left[t_0:t_1:t_0\right] \\
l_{13} &= \left[0:0:1\right] &
l_{23} &= \left[s_0 t_0:s_0 t_1:s_1 t_0\right] &
l_{33} &= \left[s_0 t_1:s_0 t_1:s_1 t_0\right]
\end{align*}
This explicit description shows that the family of $(3,3)$-nets we
have constructed is two-dimensional, parameterized by 
$\mathbb{P}^1 \times \mathbb{P}^1 = 
\{ \left[s_0 : s_1\right] \times \left[t_0 : t_1\right]\}$.
It is clear that for generic choice of $\left[s_0 : s_1\right]$ and
$\left[ t_0 : t_1 \right]$ this will be a proper $(3,3)$-net
consisting of three \emph{nondegenerate} triangles; 
and it is also clear that the parameters can be chosen
so that the configuration is in fact real.

From Yuzvinsky's construction described in the previous section, we
might have expected to find a three-dimensional family of $(3,3)$-nets,
since any such corresponds to the group $\mathbb{Z}/3\mathbb{Z}$.  What
happens here, though, is that the nine lines of the configuration do 
not impose independent conditions on smooth cubics $E$; the nine points
in the dual plane are the intersection of two cubics, and so there is
a pencil of cubics containing them.  This accounts for the ``missing''
parameter.

\subsection{$4$-nets of cubics}
Let us begin by examining something that happened in the preceding
example.  It occurred that as soon as our two triangles, defined by
$\mathcal{A}_1$ and $\mathcal{A}_2$, were perspective from both 
$l_{31}$ and $l_{32}$, they were automatically perspective from
a third line, $l_{33}$, which completed the perspectivities described
in the Latin square $M_3$.  There is a coordinate-free explanation
of this phenomenon, which we will describe now.

For simpler notation, let us write $\mathcal{A}_1 = \{ l, m, n \}$
and $\mathcal{A}_2 = \{ l', m', n' \}$.  Again we take the Latin
square
\[
M_3 = \left[ \begin{array}{ccc}
		1 & 2 & 3 \\
		2 & 3 & 1 \\
		3 & 1 & 2
		\end{array} \right] ,
\]
and (using the classical notation for perspectivities) we assume
that there exist lines $l''$ and $m''$ such that
\[
\frac{\begin{array}{ccc}
	l & m & n \\
	l' & n'& m'
	\end{array}}{l''},  \,\,\,\,\,\,\,
\frac{\begin{array}{ccc}
	l & m & n \\
	m' & l'& n'
	\end{array}}{m''} .
\]
With obvious labeling, this means we have accounted for all the 1's
and 2's appearing in $M_3$.
Now, there is a projective transformation from the pencil of points
on $l''$ to the pencil of points on $m''$ that sends
\begin{eqnarray*}
l\cap l' &\longleftrightarrow& n\cap n' \\
m\cap n' &\longleftrightarrow& l\cap m' \\
n\cap m' &\longleftrightarrow& m\cap l' .
\end{eqnarray*}
It is a theorem of elementary projective geometry that the cross joins
of such a transformation lie on a line, called the \emph{axis of 
homology}.  If we take $n''$ to be this line, we see that the three
points $l\cap n'$, $m\cap m'$, and $n\cap l'$ lie on $n''$; in other
notation, we may write
\[
\frac{\begin{array}{ccc}
	l & m & n \\
	n' & m'& l'
	\end{array}}{n''} ,
\]
which accounts for the 3's appearing in the Latin square 
$M_3$.  This explains why
$M_3$ imposes only two conditions on a pair of
triangles.

We turn now to $(4,3)$-nets.  According to 
Proposition~\ref{prop:knet-construct}, we can attempt to construct such
a configuration by starting with two orthogonal Latin squares of order
three.  Of course, we have already seen such a pair, in 
Example~\ref{ex:43net}:
\[
M_3 = \left[\begin{array}{ccc}
1 & 2 & 3 \\
2 & 3 & 1 \\
3 & 1 & 2
\end{array}\right], \,\,\,
M_4 = \left[\begin{array}{ccc}
1 & 2 & 3 \\
3 & 1 & 2 \\
2 & 3 & 1
\end{array}\right].
\]
If we start as in the preceding section, with
\begin{align*}
l_{11} &= \left[1:0:0\right] &
l_{21} &= \left[1:1:1\right] \\
l_{12} &= \left[0:1:0\right] &
l_{22} &= \left[s_0:s_1:s_2\right] \\
l_{13} &= \left[0:0:1\right] &
l_{23} &= \left[t_0:t_1:t_2\right],
\end{align*}
we see that we have four parameters to work with.  We have already seen
that $M_3$ imposes only two conditions on these parameters, and an analogous
argument shows that $M_4$ also imposes only two conditions.  Therefore we
expect that there is an essentially
unique choice of parameters $\left[s_i\right]$ and 
$\left[t_i\right]$ that will give a $(4,3)$-net.  This expectation is
correct; we omit the algebra, but it is not difficult to compute that
the coordinates of this unique $(4,3)$-net are as follows:
\begin{align*}
l_{11} &= \left[1:0:0\right] &
l_{21} &= \left[1:1:1\right] &
l_{31} &= \left[\omega:1:1\right] &
l_{41} &= \left[\omega^2:1:1\right] \\
l_{12} &= \left[0:1:0\right] &
l_{22} &= \left[1 : \omega : \omega^2 \right] &
l_{32} &= \left[1:\omega:1\right] &
l_{42} &= \left[1:\omega^2:1\right] \\
l_{13} &= \left[0:0:1\right] &
l_{23} &= \left[1 : \omega^2 : \omega \right] &
l_{33} &= \left[1:1:\omega\right] &
l_{43} &= \left[1:1:\omega^2\right] ,
\end{align*}
where $\omega$ is a primitive third root of unity.  Note that this 
configuration is not projectively equivalent to a configuration with
real coordinates.  This is not surprising: see \cite{cordovil:no-4nets} for
a proof that there are no $4$-nets realizable in the real projective plane.
(The present author thanks Professor Michael Falk for this reference.)

This configuration of lines is of course very well known: these four
(nondegenerate) triangles are the only singular members of the famous
\emph{Hesse pencil of cubics}, defined in these coordinates by
\[
\lambda xyz + \mu (x^3 + y^3 + z^3) = 0 .
\]
The four triangles correspond to
\[
\left[ \lambda : \mu \right] = \left[1:0\right], \,\,\,
\left[0:1\right], \,\,\, 
\left[6\omega + 3 : \omega\right], \,\,\,
\left[6\omega + 3 : -\omega^2\right],
\]
respectively.

We conclude this section with an observation regarding perspective
triangles.  We remarked in an earlier section that there can be at
most one perspectivity between two $d$-gons corresponding to each
permutation of the set $\{1, 2, \ldots, d\}$.  In the present case,
$d=3$, and so we see that two triangles can be perspective from at
most six lines, each corresponding to one of the six permutations
of $\{1,2,3\}$ appearing in either $M_3$ or $M_4$.
But we also observed that if two of these permutations
come from the same Latin square, then the triangles are automatically
perspective according to the third permutation in that Latin square,
as well.  This means that we have proven the following elementary
result:
\begin{prop}
Two triangles in the plane may be perspective from exactly zero, one, two,
three, four, or six points --- but not five.
\end{prop}
In fact, we have seen that there is a four-dimensional family of pairs
of triangles, and the general member of this family has no perspectivities
at all.  Imposing one perspectivity corresponds to the vanishing of one
determinant, and so defines a three-dimensional subvariety with six
components.  Imposing a
second condition from the same Latin square actually results in threefold
perspectivity, while imposing a second condition from the other Latin square
results in twofold perspectivity; these two cases correspond to different
sets of components of a two-dimensional subvariety.  
The only way to get exactly
fourfold perspectivity is to impose all conditions from one Latin square,
and just one from the other Latin square; this is a one-dimensional
subvariety with six components.  Finally, we get the essentially unique
$(4,3)$-net described above as the only example of a pair of triangles
which are in sixfold perspective.

\subsection{$3$-nets of quartics}
To construct $(3,4)$-nets, we adopt the same approach as for $(3,3)$-nets:
namely, we start with a pair of (nondegenerate) quadrangles, label their
sides, and choose a Latin square $M_3$ of order four to impose conditions
on the pair of quadrangles.  The interesting point here is that there
are two different isotopy classes of Latin squares of order four:
\[
M_3 = \left[ \begin{array}{cccc}
	1 & 2 & 3 & 4 \\
	2 & 3 & 4 & 1 \\
	3 & 4 & 1 & 2 \\
	4 & 1 & 2 & 3
	\end{array}\right], \,\,\,
M'_3 = \left[ \begin{array}{cccc}
	1 & 2 & 3 & 4 \\
	2 & 1 & 4 & 3 \\
	3 & 4 & 1 & 2 \\
	4 & 3 & 2 & 1
	\end{array}\right] .
\]
The fact that these two Latin squares represent different isotopy classes
means that there is no way to relabel the lines of a $(3,4)$-net corresponding
to $M_3$ so that it corresponds to $M'_3$, and vice versa; thus we may
expect to find two different families of $(3,4)$-nets.  This is indeed
what happens.

Up to a choice of symbols, $M_3$ is the multiplication table for
$\mathbb{Z}/4\mathbb{Z}$, and $M'_3$ is the multiplication table for
$\mathbb{Z}/2\mathbb{Z} \times \mathbb{Z}/2\mathbb{Z}$.
Both of these groups are realizable as subgroups of elliptic curves,
and so Yuzvinsky's construction applies.  We expect that the dimension
of each family is exactly three: certainly we expect at most three,
by Yuzvinsky's construction, and in fact we expect no ``discrepancy''
as in earlier examples, because there should be exactly one plane
cubic passing through the twelve points dual to a $(3,4)$-net.
Of course, this expectation is based on the assumption that any $(3,4)$-net
is algebraic, in the sense that the dual configuration of points lie on a
plane cubic, as in Yuzvinsky's construction.  
This is indeed the case, for both families of $(3,4)$-nets
(see \cite{schroeter}, for example);
but in any event, we can verify by direct computation that both families
of $(3,4)$-nets are three-dimensional, parameterized by 
$\mathbb{P}^1 \times \mathbb{P}^1 \times \mathbb{P}^1$.

Here are coordinates for a family of $(3,4)$-nets corresponding to 
the Latin square $M_3$.  Notice that for generic choice of parameters,
the three quadrangles are nondegenerate, and notice also that the 
parameters can be chosen so that the configuration is real:
\begin{align*}
l_{11} &= \left[1:1:1\right] &
l_{21} &= \left[1:0:0\right] &
l_{31} &= \left[s_0:s_1:s_1\right] \\
l_{12} &= \left[t_1 u_1:t_0 u_0:t_1 u_0\right] &
l_{22} &= \left[0:1:0\right] &
l_{32} &= \left[t_1:t_0:t_1\right] \\
l_{13} &= \left[s_0 u_1:s_1 u_1:s_1 u_0\right] &
l_{23} &= \left[0:0:1\right] &
l_{33} &= \left[u_1:u_1:u_0\right] \\
l_{14} &= \left[s_0 t_1:s_0 t_0:s_1 t_1\right] &
l_{24} &= \left[x_0:x_1:x_2\right] &
l_{34} &= \left[s_0 t_1 u_1:s_0 t_0 u_0:s_1 t_1 u_0\right] ,
\end{align*}
where we have
\begin{eqnarray*}
x_0 &=& t_1(s_0 t_1 u_1 + s_1 t_1 u_1 - s_0 t_0 u_1 - s_0 t_0 u_0), \\
x_1 &=& s_0 t_0^2 u_0 - s_1 t_1^2 u_1, \\
x_2 &=& t_1(s_1 t_1 u_0 + s_1 t_1 u_1 - s_1 t_0 u_0 - s_0 t_0 u_0) .
\end{eqnarray*}
We remark that this is a much-studied configuration, although perhaps
not well-known in modern times.  Schroeter \cite{schroeter} discusses
in great detail the geometry of this configuration, which he refers to 
as the $(12_4, 16_3)B$ configuration; the ``$B$'' is to distinguish it 
from Hesse's $(12_4, 16_3)$
configuration, which he examines in detail in the same paper.

Here are coordinates for a family of $(3,4)$-nets corresponding to 
$M'_3$.  Again, the quadrangles can be chosen to be nondegenerate and
real:
\begin{align*}
l_{11} &= \left[1:1:1\right] &
l_{21} &= \left[1:0:0\right] &
l_{31} &= \left[s_0:s_1:s_1\right] \\
l_{12} &= \left[s_0 t_1:s_1 t_0:s_1 t_1\right] &
l_{22} &= \left[0:1:0\right] &
l_{32} &= \left[t_1:t_0:t_1\right] \\
l_{13} &= \left[s_0 u_1:s_1 u_1:s_1 u_0\right] &
l_{23} &= \left[0:0:1\right] &
l_{33} &= \left[u_1:u_1:u_0\right] \\
l_{14} &= \left[t_1 u_1:t_0 u_1:t_1 u_0\right] &
l_{24} &= \left[x_0:x_1:x_2\right] &
l_{34} &= \left[s_0 t_1 u_1:s_1 t_0 u_1:s_1 t_1 u_0\right] ,
\end{align*}
where now we have
\begin{eqnarray*}
x_0 &=& (s_0 + s_1)t_1 u_1, \\
x_1 &=& s_1 (t_0 + t_1) u_1, \\
x_2 &=& s_1 t_1 (u_0 + u_1) .
\end{eqnarray*}
This is a $(12_4, 16_3)$ configuration which was once referred to as the 
\emph{Hesse configuration}, although these days \cite{dolgachev} 
it is sometimes 
referred to as the \emph{Hesse-Salmon configuration}, perhaps to 
distinguish it from the $(4,3)$-net arising from the Hesse pencil of
cubics.  The geometry of this $(12_4, 16_3)$ 
configuration is studied in detail in
both \cite{schroeter} and \cite{mathews}.  The latter reference 
is especially interesting, because Mathews relates this configuration to the
geometry of so-called \emph{desmic surfaces} in $\mathbb{P}^3$; we will
briefly describe one connection here.

Let $X, Y, Z$ be the three tetrahedra in $\mathbb{P}^3$ whose vertices
are as follows:
\begin{align*}
X_1 &= \left[1:0:0:0\right] &
Y_1 &= \left[1:1:1:1\right] &
Z_1 &= \left[-1:1:1:1\right] \\
X_2 &= \left[0:1:0:0\right] &
Y_2 &= \left[1:1:-1:-1\right] &
Z_2 &= \left[1:-1:1:1\right] \\
X_3 &= \left[0:0:1:0\right] &
Y_3 &= \left[1:-1:1:-1\right] &
Z_3 &= \left[1:1:-1:1\right] \\
X_4 &= \left[0:0:0:1\right] &
Y_4 &= \left[1:-1:-1:1\right] &
Z_4 &= \left[1:1:1:-1\right]
\end{align*}
These are called \emph{desmic tetrahedra}, and they 
satisfy a number of geometric properties in relation to each other;
the most notable of these properties --- and the one normally taken as
the definition of ``desmic'' --- is that any two are in perspective
from each vertex of the third.  These perspectivities are encoded
in the Latin square $M'_3$, as is easily checked.  It turns out that
these three tetrahedra are completely reducible members of a pencil
of quartic surfaces in $\mathbb{P}^3$, called \emph{desmic surfaces};
a generic hyperplane section will give a pencil of quartics with 
three completely reducible members; these three reducible members
form a $(3,4)$-net corresponding to $M'_3$.  This gives another way
of seeing that this family of $(3,4)$-nets is three-dimensional: any
three desmic tetrahedra are projectively equivalent to the one given
above in coordinates, and then the choice of a hyperplane section gives
us the three parameters.

\subsection{$3$-nets of quintics}
Let us start by considering the Latin square which, up to choice of 
symbols, is the multiplication table of the group $\mathbb{Z}/5\mathbb{Z}$:
\[
M_3 = \left[ \begin{array}{ccccc}
	1 & 2 & 3 & 4 & 5 \\
	2 & 3 & 4 & 5 & 1 \\
	3 & 4 & 5 & 1 & 2 \\
	4 & 5 & 1 & 2 & 3 \\
	5 & 1 & 2 & 3 & 4
	\end{array} \right] .
\]
We can carry out the same sort of analysis as in previous examples to 
find a family of $(3,5)$-nets, consisting of three nondegenerate 
pentagons; with the aid of a computer algebra package, we can carry
out the elimination theory to get the following coordinates for the
lines of the configuration:
\begin{align*}
l_{11} &= \left[1:1:1\right] &
l_{21} &= \left[1:0:0\right] &
l_{31} &= \left[s_2 t_0:s_1 t_2:s_1 t_2\right] \\
l_{12} &= \left[s_0 t_1:s_1 t_1:s_2 t_0\right] &
l_{22} &= \left[0:1:0\right] &
l_{32} &= \left[s_2 t_0:s_1 t_1:s_2 t_0\right] \\
l_{13} &= \left[s_0 t_0:s_0 t_1:s_2 u_0\right] &
l_{23} &= \left[0:0:1\right] &
l_{33} &= \left[s_0 t_1:s_0 t_1:s_2 t_0\right] \\
l_{14} &= \left[s_2 t_0:s_1 t_2:s_2 t_2\right] &
l_{24} &= \left[x_0:x_1:x_2\right] &
l_{34} &= \left[s_0:s_1:s_2\right] \\
l_{15} &= \left[s_2 t_0:s_1 t_1:s_1 t_2\right] &
l_{25} &= \left[y_0:y_1:y_2\right] &
l_{35} &= \left[t_0:t_1:t_2\right] 
\end{align*}
Here the parameters $\left[s_i\right]$ and $\left[t_i\right]$ are not free;
they must be taken on the degree $(3,3)$ hypersurface in
$\mathbb{P}^2 \times \mathbb{P}^2$ defined by the single equation
\[
\begin{split}
s_0^2 s_1 t_1^2 t_2 
&-s_0^2 s_1 t_1 t_2^2 - s_0 s_1^2 t_0 t_1 t_2 + s_0 s_1^2 t_1 t_2^2 \\
&+s_0 s_1 s_2 t_0^2 t_2 - s_0 s_1 s_2 t_1^2 t_2 - s_0 s_2^2 t_0^2 t_1 \\
&+s_0 s_2^2 t_0 t_1 t_2 + s_1 s_2^2 t_0^2 t_1 - s_1 s_2^2 t_0^2 t_2 = 0 .
\end{split}
\]
The expressions for $\left[x_i\right]$ and $\left[y_i\right]$ in terms
of $\left[s_i\right]$ and $\left[t_i\right]$ are politely omitted; they
are easily calculated from the relevant perspectivities:
\[
\frac{\begin{array}{ccccc}
	l_{11} & l_{12} & l_{13} & l_{14} & l_{15} \\
	l_{34} & l_{35} & l_{31} & l_{32} & l_{33}
	\end{array}}{l_{24}},  \,\,\,\,\,\,\,
\frac{\begin{array}{ccccc}
	l_{11} & l_{12} & l_{13} & l_{14} & l_{15} \\
	l_{35} & l_{31} & l_{32} & l_{33} & l_{34}
	\end{array}}{l_{25}} .
\]
So once again we have a three-dimensional family of configurations.
This is the dimension we would expect from Yuzvinsky's construction if
we knew \emph{a priori} that any $(3,5)$-net corresponding to $M_3$ is
algebraic; the present author does not know a direct proof of this fact,
however.

We note that for every $3$-net we have considered, its corresponding 
Latin square is isotopic to the multiplication table of a group.
This phenomenon led Yuzvinsky to ask whether there are any $3$-nets 
whose Latin square is \emph{not} isotopic to the multiplication table
of a group \cite{yuz:nets}.  The following example, which the present
author believes to be new, answers Yuzvinsky's question affirmatively.

We begin with the following Latin square:
\[
M'_3 = \left[ \begin{array}{ccccc}
	1 & 2 & 3 & 4 & 5 \\
	2 & 1 & 4 & 5 & 3 \\
	3 & 5 & 1 & 2 & 4 \\
	4 & 3 & 5 & 1 & 2 \\
	5 & 4 & 2 & 3 & 1
	\end{array} \right] .
\]
This Latin square is not isotopic to one which is a multiplication
table of a group; in fact, $M_3$ and $M'_3$ represent the only two
isotopy classes of Latin squares of order five.  As before, we carry
out straightforward but tedious analysis to find the following
coordinates for the lines of a $(3,5)$-net corresponding to $M'_3$:
\begin{align*}
l_{11} &= \left[1:1:1\right] &
l_{21} &= \left[1:0:0\right] &
l_{31} &= \left[s_2 t_1:s_1 t_2:s_1 t_2\right] \\
l_{12} &= \left[s_0 s_2 t_1:s_1 s_2 t_1: s_0 s_1 t_2\right] &
l_{22} &= \left[0:1:0\right] &
l_{32} &= \left[s_0 t_2:s_2 t_1:s_0 t_2\right] \\
l_{13} &= \left[s_2 t_0 t_1:s_1 t_0 t_2:s_2 t_1 t_2\right] &
l_{23} &= \left[0:0:1\right] &
l_{33} &= \left[s_1 t_0:s_1 t_0:s_2 t_1\right] \\
l_{14} &= \left[s_1 t_0:s_1 t_1:s_2 t_1\right] &
l_{24} &= \left[x_0:x_1:x_2\right] &
l_{34} &= \left[s_0:s_1:s_2\right] \\
l_{15} &= \left[s_0 t_2:s_2 t_1:s_2 t_2\right] &
l_{25} &= \left[y_0:y_1:y_2\right] &
l_{35} &= \left[t_0:t_1:t_2\right] 
\end{align*}
Again, the parameters $\left[s_i\right]$ and $\left[t_i\right]$ are
not free; they must be chosen to lie on the degree $(3,3)$ hypersurface
in $\mathbb{P}^2 \times \mathbb{P}^2$ defined by the equation
\[
s_0 s_1^2 t_0 t_2^2 - s_0 s_1 s_2 t_1 t_2^2 - s_1^2 s_2 t_0 t_1 t_2
+s_1^2 s_2 t_1 t_2^2 - s_1 s_2^2 t_1^2 t_2 + s_2^3 t_1^3 = 0 .
\]
Again we omit the explicit formulas for $\left[x_i\right]$ and 
$\left[y_i\right]$; they may be computed from the 
following perspectivities:
\[
\frac{\begin{array}{ccccc}
	l_{11} & l_{12} & l_{13} & l_{14} & l_{15} \\
	l_{34} & l_{35} & l_{32} & l_{31} & l_{33}
	\end{array}}{l_{24}},  \,\,\,\,\,\,\,
\frac{\begin{array}{ccccc}
	l_{11} & l_{12} & l_{13} & l_{14} & l_{15} \\
	l_{35} & l_{33} & l_{34} & l_{32} & l_{31}
	\end{array}}{l_{25}} .
\]
And so we have a three-dimensional family of configurations
corresponding to the Latin square $M'_3$.

\section{Conclusion}
Since there seems to be no reason to assume that the Latin square
corresponding to a $3$-net is isotopic to the multiplication
table of a group, it would be interesting to know whether there
are families of $3$-nets corresponding to \emph{every} isotopy
class of Latin squares.  This question is not easy to dismiss
by parameter count alone.  The phenomenon we discussed in the 
sections on $(3,3)$ and $(4,3)$-nets occurs in general: namely,
the various perspectivities imposed by a Latin square are not
independent, as may be checked in our other examples.

On the other hand, the present author is not bold enough to 
conjecture that every isotopy class of Latin squares gives a
corresponding family of $3$-nets in the complex projective
plane.  Without any deep insight into the question, the only
technique available is brute force elimination theory --- and
already in degree six, the force required is too brute for a
computer algebra package.  So we will close
this paper with a question of our own:

\begin{prob}
Are there isotopy classes of Latin squares of order~$d \ge 6$ 
which are not realized by a three-dimensional 
family of $(3,d)$-nets in the complex projective plane?
\end{prob}

\end{document}